\documentclass[12pt, paperwidth=5.35in, paperheight=8.3in]{amsart}

\pdfpagewidth=\paperwidth
\pdfpageheight=\paperheight

\newtheorem{theorem}{Theorem}[section]
\newtheorem{lemma}[theorem]{Lemma}

\theoremstyle{definition}
\newtheorem{definition}[theorem]{Definition}

\theoremstyle{remark}
\newtheorem{remark}[theorem]{Remark}
\newtheorem{corollary}[theorem]{Corollary}
\newtheorem{proposition}[theorem]{Proposition}
\numberwithin{equation}{section}

\allowdisplaybreaks
\usepackage{amsmath,amssymb,amscd,amsthm,amsxtra}
\usepackage[colorlinks=true, pdfstartview=FitV, linkcolor=blue, citecolor=blue, urlcolor=blue]{hyperref}
\usepackage{color}
\usepackage{graphicx,epstopdf,verbatim}

\newcommand{\q}{\quad}

\newcommand{\R}{\mathbb{R}}

\newcommand{\ra}{\rightarrow}

\begin{document}
\title[MULTIVARIABLE DYADIC PARAPRODUCTS]{\textbf{WEIGHTED INEQUALITIES FOR MULTIVARIABLE DYADIC PARAPRODUCTS }}
\author[D. CHUNG]{DAEWON CHUNG}
\address{Department of mathematics and Statistics, 1 University of New Mexico, Albuquerque, NM 87131-001}
\email{dwchung@unm.edu}
\date{\today}
 \maketitle
\begin{abstract}
Using Wilson's Haar basis in $\R^n$, which is different than the usual tensor product Haar functions,
we define its associated dyadic paraproduct in $\R^n$. We can then extend ``trivially'' Beznosova's Bellman function proof of the
linear bound in $L^2(w)$ with respect to $[w]_{A_2}$ for the 1-dimensional dyadic paraproduct. Here trivial means that each piece of the argument
that had a Bellman function proof has an $n$-dimensional counterpart that holds with the same Bellman function. The lemma that allows for this painless
extension we call the good Bellman function Lemma. Furthermore the argument allows to obtain dimensionless bounds in the anisotropic case.
\end{abstract}

\maketitle

\section{Introduction and Main results}
The name \emph{Paraproduct} was coined by Bony, in 1981 (see
\cite{Bo}), who used paraproducts to linearize the problem in the
study of singularities of solutions of semilinear partial
differential equations. After his work, the paraproducts have
played an important role in harmonic analysis because they are
examples of singular integral operators which are not
translation-invariant. Also, every
singular integral operator which is bounded on $L^2$ decomposes
into a paraproduct, an adjoint of a paraproduct, and an almost convolution operator. Moreover they arise
as building blocks for more general operators such as multipliers.

For the locally integrable functions $b$ and $f$, the dyadic
paraproduct is defined by
$$\pi_bf:=\sum_{I\in\mathcal{D}}\langle b,h_I\rangle\langle
f\rangle_I h_I\,,$$ on the real line. Here the $\mathcal{D}$ denotes the collection of all dyadic intervals.
$\{h_I\}_{I\in\mathcal{D}}$ is the Haar basis in $L^2_{\R}$, $\langle\cdot,\cdot\rangle$ stands for the standard inner product in $L^2_{\R}$, and
$\langle \cdot\rangle$ denotes the average over the interval $I$. It is now well known fact that the dyadic
paraproduct is bounded on $L^p$ if $b\in BMO^d$ (see \cite{Per1}).
Thus, after we fix $b$ in $BMO^d$, we consider $\pi_bf$ as a
linear operator acting on $f$. We say the positive almost everywhere and locally integrable function $w$, a weight, satisfies the $A_p$ condition
if: \begin{equation}[w]_{A_p}:=\sup_{I}\langle w\rangle_I\langle w^{-1/(p-1)}\rangle^{p-1}_I<\infty\,,\label{AP}\end{equation}
where the supremum is taken over all intervals.
The class of weight $A_p$ was first
presented in \cite{MB}, the Hardy-Littlewood maximal operator is
bounded on $L^p(w)$ for $1<p<\infty$ if and only if the weight $w$
belong to the class $A_p$. If we take the supremum over all dyadic intervals in (\ref{AP}),
then we call it is the dyadic $A_p$-characteristic and is denoted by $[w]_{A_p^d}\,.$
Beznosova proved in \cite{Be} that the following linear estimate holds for the dyadic paraproduct in
$L^2(w)\,.$

\begin{theorem}[O. Beznosova]\label{lbp}
The norm of dyadic paraproduct on the weighted Lebesgue space
$L^2_{\R}(w)$ is bounded from above by a constant multiple of the
product of the $A^d_2$-characteristic of the weight $w$ and the
$BMO^{\,d}$ norm of $b$, that is
$$\|\pi_bf\|_{L^2(w)}\leq C[w]_{A^d_2}\|b\|_{BMO^d}\|f\|_{L^2(w)}\,.$$
\end{theorem}

In fact, the linear bound in $L^2_{\R^n}(w)\,$ of the $n$-dimensional dyadic paraproducts are
recovered in \cite{HyLReVa, CrMP} using different methods.  However, in this paper, to prove the
linear bound of the dyadic paraproduct in $L^2_{\R^n}(w)\,$ we use
the Bellman function arguments as in \cite{Be}.

One of the
main purposes in this paper is an estimate for the $n$-dimensional
analog of the dyadic paraproduct and to establish a linear bound
with $[w]_{A^d_2}$ and $\|b\|_{BMO^d}\,.$ Thus, throughout the
paper, we will be concerned with a class of weights, $A^d_p\,,$ on
$\R^n.$ If
$$ [w]_{A^d_p}:=\sup_{Q\in\mathcal{D}^n}\langle w\rangle_Q\langle w^{-1/(p-1)}\rangle^{p-1}_{Q}<\infty\,,$$
then we say the weight $w$ belongs to the class of $A^d_p$
weights. Here $\mathcal{D}^n$ denotes the collection of all dyadic
cubes in $\R^n.$

In order to extend Theorem
\ref{lbp} to the multivariable setting in the spirit of \cite{Be},
we are using Bellman function arguments. This allows to establish
the dimension free estimates in terms of anisotropic weight
characteristic. Thus we need to consider the class of anisotropic
$A_2$-weights and the class of anisotropic $BMO$ functions which
are defined as follows.
\begin{definition}
A locally integrable and positive almost everywhere function $w$
on the space $\R^n$ belongs to class of $A^R_p$ weights, $
1<p<\infty$ if $$[w]_{A_p^R}:=\sup_{R}\langle w\rangle_R\langle
w^{-1/(p-1)}\rangle_R^{p-1}<\infty\,,$$ where the supremum is
taken over all rectangles $R\subset\R^n$ with sides parallel to
the coordinate axes.
\end{definition}
\begin{definition}
A locally integrable function on $\R^n$ belongs to $BMO^R$ if
$$\|b\|_{BMO^R}:=\sup_{R}\frac{1}{|R|}\int_R|\,b(x)-\langle
b\rangle_R|\,dx<\infty\,,$$ where the supremum runs over all
rectangles $R\subset\R^n$ with sides parallel to the coordinate
axes.
\end{definition}
Since a cube is a particular case of a rectangle, it is easy to
observe that $\|b\|_{BMO}\leq \|b\|_{BMO^R}$.  In \cite{K}, one can find the example which is in $BMO_{\R^2}$ but not in
$BMO^R_{\R^2}\,.$ Thus, $BMO\supset
BMO^R$ when $n\geq2\,.$ It is a well known fact that the weight $|x|^{\alpha}\in A_2$ if and only if $|\alpha|<n\,.$ Thus, $|x|\in A_2$ in $\R^2\,.$ However, one can see that $\langle |x|\rangle_{R_t}\langle |x|^{-1}\rangle_{R_t}$
 behaves like $\log t$, where $R_t=[0,t]\times [0,1]\,.$ Even if $[w]_{A_2}=[w]_{A^R_2}$ in $\R$,
 we can see that the $A_2^R$ weight class belongs strictly to $A_2$, when $n\geq 2\,.$ We now state our main
results. Here $\pi_b$ is the dyadic paraproduct associated to Wilson's Haar basis in $\R^n$ and defined in Section 2.
\begin{theorem}\label{MA}
For $1<p<\infty$ there exists constants $C(n,p)$ only depending
on $p$ and dimension $n$ and $C$ which doesn't depend on the
dimensional constant such that
$$\|\pi_b\|_{L^p_{\R^n}(w)\ra L^p_{\R^n}(w)}\leq C(n,p)[w]_{A^d_p}^{\max\{1,\frac{1}{p-1}\}}\|b\|_{BMO^d_{\R^n}}\,$$
and for all weights $w\in A^d_p$
and $b\in BMO^d_{\R^n}\,,$ and
$$\|\pi_b\|_{L^2_{\R^n}(w)\ra
L^2_{\R^n}(w)}\leq C\,[w]_{A^R_2}\|b\|_{BMO^R_{\R^n}}\,$$ for all weight $w\in A^R_2$ and
$b\in BMO^R_{\R^n}\,.$
\end{theorem}

Through out the paper, we denote a constant by $C$ which may
change line by line and we indicate its dependence on parameters using a
parenthesis, for example $C(n,p)\,.$ In Section 2 we will discuss some $n$-dimensional Haar systems and
introduce notations. In Section 3 we will introduce the multivariable dyadic paraproducts. In Section 4 we will introduce certain
embedding theorems and weight inequalities which are extended to
the several variable setting. In Section 5 we prove the main
results which provide the linear bounds for the dyadic paraproduct
in $L^2_{\R^n}(w)\,,$ and dimension free estimates. We remark in Section 6 that similar method recovers known estimates for martingale transform and
obtains dimension free estimates as well.

\section{Wilson's Haar system in $\R^n$}
First of all, we need to introduce the appropriate $n$-dimensional Haar
systems. For any $Q\in\mathcal{D}^n\,,$ we set
$\mathcal{D}_1^n(Q)\equiv\{Q'\in\mathcal{D}^n:\,Q'\subset Q,\,
\ell(Q')=\ell(Q)/2\}\,,$ the class of $2^n$ dyadic sub-cubes of
$Q$\,, where we denote the side length of cubes by $\ell(Q)\,.$ We will also denote the class of all dyadic
sub-cubes of $Q$ by $\mathcal{D}^n(Q)$. Then we can write $\mathcal{D}^n(Q)=\bigcup_{j=0}^{\infty}\mathcal{D}_j^n(Q)\,.$ We
refer to \cite{Wil} for the following lemma.

\begin{lemma}\label{PDC}
Let $Q\in\mathcal{D}^n$. Then, there are $2^n-1$ pairs of sets
$\{(E^1_{j,Q},\,E^2_{j,Q})\}_{j=1}^{2^n-1}$ such that:
\begin{enumerate}
\item[(1)] for each $j\,,$
$\big|E^1_{j,Q}\big|=\big|E^2_{j,Q}\big|\,.$ \item[(2)] for each
$j\,,$ $E^1_{j,Q}$ and $E^2_{j,Q}$ are non-empty unions of cubes from
$\mathcal{D}_1^n(Q)$; \item[(3)] for each $j\,,$ $E^1_{j,Q}\cap
E^2_{j,Q}=\emptyset$; \item[(4)] for every $j\neq k,$ one of the
following must hold:
\begin{enumerate}
\item[(a)] $E^1_{j,Q}\cup E^2_{j,Q}$ is entirely contained in
either $E^1_{k,Q}$ or $E^2_{k,Q}$; \item[(b)] $E^1_{k,Q}\cup
E^2_{k,Q}$ is entirely contained in either $E^1_{j,Q}$ or
$E^2_{j,Q}$; \item[(c)] $\big(E^1_{j,Q}\cup E^2_{j,Q}\big) \cap
\big(E^1_{k,Q}\cup E^2_{k,Q}\big)=\emptyset\,.$
\end{enumerate}
\end{enumerate}
\end{lemma}

We can construct such a set by induction on $n\,.$ It is clear
when $n=1\,.$ We assume that Lemma \ref{PDC} is true for $n-1$ and
let $\widetilde{Q}$ be the $(n-1)$-dimensional cube and
$\big\{(E^1_{j,\widetilde{Q}},E^2_{j,\widetilde{Q}})\big\}_{j=1}^{2^{n-1}-1}$
be the corresponding pairs of sets for $\widetilde{Q}\,.$ We can
get the first pair of sets by $(E^1_{1,Q},E^2_{1,Q}):=(\widetilde{Q}\times
I_-,\widetilde{Q}\times I_+)$ where $I$ is a dyadic interval so
that $|I|=\ell(\widetilde{Q})\,,$ and $\widetilde{Q}\times I=Q\,.$
We also have the last $2^n-2$
pairs of sets as follows.
\begin{align*}
 \{(E^1_{2j,Q},&E_{2j,Q}^2),(E_{2j+1,Q}^1,E_{2j+1,Q}^2)\}_{j=1}^{2^{n-1}-1}\\
&:=\big\{(E^1_{j,\widetilde{Q}}\times I_-,E^2_{j,\widetilde{Q}}\times I_-),\,(E^1_{j,\widetilde{Q}}\times I_+,E^2_{j,\widetilde{Q}}\times I_+)\big\}_{j=1}^{2^{n-1}-1}\,.
\end{align*}
To save space, we denote $E^1_{j,Q}\cup E^2_{j,Q}$ by $E_{j,Q}\,$ and, by $(1)$ in Lemma \ref{PDC}, we have $|E_{j,Q}|=2|E^i_{j,Q}|$ for $i=1,2\,.$
Note that the sets $E_{j,Q}$ are rectangles. Also note that we assign $E_{1,Q}=Q\,,$ $E_{2,Q}=E^1_{1,Q}$ and $E_{3,Q}=E^2_{1,Q}$ and so on.
With such a choice, we have
\begin{align*}
Q&=E_{1,Q}=E_{1,Q}^1\cup E_{1,Q}^2=E_{2,Q}\cup E_{3,Q}=(E_{4,Q}\cup E_{6,Q})\cup (E_{5,Q}\cup\ E_{7,Q})\\
&=\cdots=E_{2^{n-1},Q}\cup E_{2^{n-1}+1,Q}\cup\cdots\cup E_{2^{n}-1,Q}\\
&=E^1_{2^{n-1},Q}\cup E^2_{2^{n-1},Q}\cup E^1_{2^{n-1}+1,Q}\cup E^2_{2^{n-1}+1,Q}\cup \cdots\cup E^1_{2^{n}-1,Q}\cup E^2_{2^{n}-1,Q}\,,
\end{align*}
in fact, $$Q=\bigcup_{j=2^k}^{2^k-1}E_{j,Q}=\bigcup_{j=2^k}^{2^{k}-1}(E^1_{j,Q}\cup E^2_{j,Q})\,,$$
the sets $E_{j,Q}$ in that range of $j$'s are disjoint, and $$\mathcal{D}^n_1(Q)=\{E^1_{2^{n-1},Q}\,, E^2_{2^{n-1},Q}\,, E^1_{2^{n-1}+1,Q}\,, E^2_{2^{n-1}+1,Q}\,, ...\,, E^1_{2^{n}-1,Q}\,, E^2_{2^{n}-1,Q}\}\,.$$

As a consequence of Lemma \ref{PDC}, we can introduce the proper
weighted Wilson's Haar system for $L^2_{\R^n}(w)\,,$ $\{h^w_{j,Q}\}_{1\leq
j\leq 2^n-1,Q\in\mathcal{D}^n}\,,$ where
 $$h^w_{j,Q}:=\frac{1}{\sqrt{w\big(E_{j,Q}\big)}}\Bigg[\sqrt{\frac{w\big(E^1_{j,Q}\big)}{w\big(E^2_{j,Q}\big)}}\chi_{E^2_{j,Q}}-\sqrt{\frac{w\big(E^2_{j,Q}\big)}{w\big(E^1_{j,Q}\big)}}\chi_{E^1_{j,Q}}\Bigg]\,.$$
When $w\equiv 1$ we denote the Wilson's Haar functions by $h_{j,Q}$. Then, every function $f\in L^2_{\R^n}(w)$ can be written as $$f=\sum_{Q\in\mathcal{D}^n}\sum_{j=1}^{2^n-1}\langle f,h^w_{j,Q}\rangle_wh^w_{j,Q}\,.$$
Moreover, $\|f\|^2_{L^2_{\R^n}(w)}=\sum_{Q\in\mathcal{D}^n}\sum_{j=1}^{2^n-1}|\langle
f,h^w_{j,Q}\rangle_w|^2\,.$ For all $Q'\in\mathcal{D}_1^n(Q)\,,$ the $h_{j,Q}$'s and $h^w_{j,Q}$'s
are constant on $Q'\,,$ we will also denote this constant by $h_{j,Q}(Q')$ and
$h^w_{j,Q}(Q')$ respectively. We
can obtain the weighted average of $f$ over $E_{j,Q}$ for some $1\leq j\leq 2^n-1,$
\begin{equation}
\langle f\rangle_{E_{j,Q},w}=\sum_{Q'\in\mathcal{D}^n:Q'\supseteq Q}\sum_{i:E_{i,Q'}\supsetneq E_{j,Q}}\langle f,h^w_{i,Q'}\rangle h^w_{i,Q'}(E_{j,Q})\,.\label{Pr:A}
\end{equation}
Furthermore, for $j=1$, $E_{1,Q}=Q,$ we have
\begin{equation}
\langle f\rangle_{E_{1,Q},w}=\langle f\rangle_{Q,w}=\sum_{Q'\in\mathcal{D}^n:Q'\supsetneq
Q}\sum_{j=1}^{2^n-1}\langle
f,h^w_{j,Q'}\rangle_wh^w_{j,Q'}(Q)\,.\label{PN:e4}
\end{equation}
Because it is occasionally more convenient to deal with simpler
functions, it might be good to have an orthogonal system in
$L^2_{\R^n}(w)$. Let us
define
\begin{equation}H^w_{j,Q}:=h_{j,Q}\sqrt{\big|E_{j,Q}\big|}-A^w_{j,Q}\,\chi_{E_{j,Q}}\,,
\textrm{ where }A^w_{j,Q}:=\frac{\langle w\rangle_{E^2_{j,Q}}-\langle w\rangle_{E^1_{j,Q}}}{2\langle w\rangle_{E_{j,Q}}}\,.\label{PN:e5}\end{equation}
Then, the family of functions $\{w^{1/2}H^w_{j,Q}\}_{j,Q}$ is an
orthogonal system for $L^2_{\R^n}$ with norms satisfying the
inequality $$\|w^{1/2}H^w_{j,Q}\|_{L^2_{\R^n}}\leq
\sqrt{\big|E_{j,Q}\big|\langle w\rangle_{E_{j,Q}}}\,.$$ By Bessel's inequality in $L^2_{\R^n}$ one gets, for all $g\in
L^2_{\R^n}\,,$
\begin{equation}\sum_{Q\in\mathcal{D}^n}\sum_{j=1}^{2^n-1}\frac{\langle gw^{1/2},H^w_{j,Q}\rangle^2}{|E_{j,Q}|\,\langle w\rangle_{E_{j,Q}}}\leq \|g\|^2_{L^2_{\R^n}}\,.\label{PN:e7}\end{equation}

As well as in the one dimensional case, one can define
\begin{equation}\|b\|_{BMO^d_{\R^n}}:=\sup_{Q\in\mathcal{D}^n}\frac{1}{|Q|}\int_Q|b(x)-\langle b\rangle_Q|dx\,,\label{BMO:e1}\end{equation}
for a locally integrable function on $\R^n\,.$ The function $b$ is
said to have dyadic bounded mean oscillation if
$\|b\|_{BMO^d_{\R^n}}<\infty\,,$ and we denote the class of all
locally integrable functions $b$ on $\R^n$ with dyadic bounded
mean oscillation by $BMO^d_{\R^n}\,.$ Notably one can replace
(\ref{BMO:e1}) by
\begin{equation}
\|b\|^2_{BMO^d_{\R^n}}=\sup_{Q\in\mathcal{D}^n}\frac{1}{|Q|}\sum_{Q\in\mathcal{D}^n(Q)}\sum_{j=1}^{2^n-1}|\langle
b,h_{j,Q}\rangle|^2\,.\label{BMO:e2}
\end{equation}
In the anisotropic case, it is known that the John-Nirenberg
inequality holds for all $b\in BMO^R$ and any rectangle
$R\subset\R^n\,$ (see \cite{K}), \begin{equation}|\{x\in R\,|\,|\,b(x)-\langle
b\rangle_R|>\lambda\}|\leq
e^{1+2/e}|R|\exp\bigg(-\frac{2/e}{\|b\|_{BMO^R}}\lambda\bigg),\quad\lambda>0\,.\label{BMOR:e1}\end{equation}
Note that the John-Nirenberg inequality is dimensionless in the
anisotropic case. As an easy consequence of (\ref{BMOR:e1}), we
have a self improving property for the anisotropic $BMO$ class.
For any rectangle $R\in\R^n\,,$ there exists a constant $C(p)$
independent of the dimension $n\,$ such that
\begin{equation}
\bigg(\frac{1}{|R|}\int_R|\,b(x)-\langle
b\rangle_R|^pdx\bigg)^{1/p}\leq
C(p)\,\|b\|_{BMO^R}\,.\label{BMOR:e2}
\end{equation}
Using the self improving property (\ref{BMOR:e2}), we have, for $i=1,...,2^n-1$, $Q'\in\mathcal{D}^n$
\begin{equation}
\frac{1}{|E_{i,Q'}|}\sum_{Q\in\mathcal{D}^n(Q')}\sum_{j:E_{j,Q}\subseteq
E_{i,Q'}}\langle b,h_{j,Q}\rangle^2\leq
C\|b\|_{BMO^R}^2\,.\label{BMOR:e3}
\end{equation}
\section{The multivariable dyadic paraproduct associated with Wilson's Haar system}
We now define the multivariable dyadic paraproduct. It is well known fact that the product of two square integrable
functions can be written as the sum of two dyadic paraproducts and
a diagonal term in a single variable case. Moreover, the diagonal
term is the adjoint of one dyadic paraproduct i.e. for all $f,g\in
L^2_{\R}\,,$
\begin{equation}fg=\pi^{\ast}_g(f)+\pi_g(f)+\pi_f(g)\,.\label{decomp}\end{equation}
Thus, we expect to have analogous decomposition. Let us assume
that $f,g\in L^2_{\R^n}\,.$ Expanding $f$ and $g$ in Wilson's Haar
system,
$$f=\sum_{Q\in\mathcal{D}^n}\sum_{j=1}^{2^n-1}\langle f,h_{j,Q}\rangle h_{j,Q}\,,\quad g=\sum_{Q'\in\mathcal{D}^n}\sum_{i=1}^{2^n-1}\langle g,h_{i,Q'}\rangle h_{i,Q'}$$
and multiplying these sums formally we can get
$$fg=\sum_{Q\in\mathcal{D}^n}\sum_{j=1}^{2^n-1}\sum_{Q'\in\mathcal{D}^n}\sum_{i=1}^{2^n-1}\langle f,h_{j,Q}\rangle\langle g,h_{i,Q'}\rangle h_{j,Q}(x)h_{i,Q'}=(I)+(II)+(III)\,.$$
Here, $(I)$ is the diagonal term $Q'=Q,~j=i$;
\begin{equation}(I)=\sum_{Q\in\mathcal{D}^n}\sum_{j=1}^{2^n-1}\langle f,h_{j,Q}\rangle\langle g,h_{j,Q}\rangle h_{j,Q}^2=\sum_{Q\in\mathcal{D}^n}\sum_{j=1}^{2^n-1}\langle f,h_{j,Q}\rangle\langle g,h_{j,Q}\rangle\frac{\chi_{E_{j,Q}}}{\big|E_{j,Q}\big|}\,.\label{PN:e10}\end{equation}
The second term $(II)$ is the upper triangle term corresponding to those $Q'\supsetneq Q$, all $i,j$ and $Q'=Q$ so that $E_{i,Q'}\supsetneq E_{j,Q}$.
\begin{equation}
(II)=\sum_{Q\in\mathcal{D}^n}\sum_{j=1}^{2^n-1}\langle f,h_{j,Q}\rangle \langle
g\rangle_{E_{j,Q}}h_{j,Q}\,,\label{PN:e11}
\end{equation}
 where we used formula (\ref{Pr:A}) for the average of $g$ on $E_{j,Q}\,.$ Similarly, the third term is the lower triangle corresponding to those $Q'\subsetneq Q$, all $i,j$ and $Q'=Q$ so that $E_{i,Q}\subsetneq E_{j,Q}\,.$
\begin{equation}
(III)=\sum_{Q\in\mathcal{D}^n}\sum_{j=1}^{2^n-1}\langle g,h_{j,Q}\rangle\langle
f\rangle_{E_{j,Q}} h_{j,Q}\,.\label{PN:e12}
\end{equation}
If we consider the sum (\ref{PN:e10}) as an operator acting on
$f,$ then we can easily check that $(III)$ is its adjoint
operator. We now can define the multivariable dyadic paraproduct by pairing
the dyadic BMO function. In $\R^n$, the dyadic paraproduct associated with Wilson's Haar system is an
operator $\pi_b\,,$ given by
\begin{equation}\pi_bf(x)=\sum_{Q\in\mathcal{D}^n}\sum_{j=1}^{2^n-1}\langle
f\rangle_{E_{j,Q}}\langle b,h_{j,Q}\rangle
h_{j,Q}(x)\,.\label{MutPa}\end{equation} Note that the
construction of the Haar systems are not unique. One can actually
construct different Haar systems \cite{DPetV}. Furthermore, the
dyadic paraproduct depends on the choice of the Haar functions.
Thus, one can establish the different dyadic paraproducts
associated with different Haar functions. But the decomposition
(\ref{decomp}) holds for all of them. We will finish this section by
including a comparison to the standard tensor product Haar basis in $\R^n$, $\{h^s_{\sigma,Q}\}$, with Wilson's Haar basis introduced in Section 2.2 and associated paraproducts. Let us denote the Haar function associated with a dyadic interval $I\in\mathcal{D}$ by $h^0_I=|\,I|^{-1/2}(\chi_{I_+}-\chi_{I_-})$ and normalized characteristic functions $h^1_I=|\,I|^{-1/2}\chi_I\,.$ Here $0$ stands for mean value zero and $1$ for the indicator. Also we consider a set of signatures $\Sigma=\{0,1\}^{\{1,...,n\}}\setminus\{(1,...,1)\}\,$ which contains $2^n-1$ signatures. These are all $n$-tuples with entries $0$ and $1$, but excluding $n$-tuple whose entries are all $1$. Then, for each dyadic cube $Q=I_1\times \cdots\times I_n$, one can get the standard tensor product Haar basis in $\R^n$ by
$$h^s_{\sigma,Q}(x_1,...,x_n)=h_{I_1}^{\sigma_1}(x_1)\times\cdots\times h_{I_n}^{\sigma_n}(x_n)\,,$$ where $\sigma=(\sigma_1,...,\sigma_n)\in\Sigma\,.$
Notice that all $h^s_{\sigma,Q}$ are supported on $Q$. In this case, we have the paraproduct associated to the standard tensor product Haar basis:
\begin{equation}\pi^s_bf(x)=\sum_{Q\in\mathcal{D}^n}\langle f\rangle_Q\sum_{\sigma\in\Sigma}\langle b,h^s_{\sigma,Q}\rangle h^s_{\sigma,Q}(x)\,.\label{MutPaS}\end{equation}
Observe that, for each dyadic cube $Q\in\mathcal{D}^n$,
\begin{equation}\mathcal{W}(Q)=\textrm{span}\{h^s_{\sigma,Q}\}_{\sigma\in\Sigma}=\textrm{span}\{h_{j,Q}\}_{j=1}^{2^n-1}\,.\end{equation}
Hence $$\textrm{Proj}_{\mathcal{W}(Q)}b=\sum_{\sigma\in\Sigma}\langle b,h^s_{\sigma,Q}\rangle h^s_{\sigma, Q}=
\sum_{j=1}^{2^n-1}\langle b,h_{j,Q}\rangle h_{j,Q}\,.$$
Changing the basis, we can see that two multivariable paraproducts, (\ref{MutPa}) and (\ref{MutPaS}), are in general different, that is
\begin{align*}
\pi^s_bf(x)&=\sum_{Q\in\mathcal{D}^n}\langle f\rangle_Q\sum_{\sigma\in\Sigma}\langle b,h^s_{\sigma,Q}\rangle h^s_{\sigma,Q}
=\sum_{Q\in\mathcal{D}^n}\langle f\rangle_Q\sum_{j=1}^{2^n-1}\langle b,h_{j,Q}\rangle h_{j,Q}\\
&\neq \sum_{Q\in\mathcal{D}^n}\sum_{j=1}^{2^n-1}\langle f\rangle_{E_{j,Q}}\langle b,h_{j,Q}\rangle h_{j,Q}=\pi_bf(x)\,.
\end{align*}

\section{Embedding theorems and weighted inequalities in $\R^n$}

In general, once we have a Bellman function proof for a certain
property in $\R$ then we can extend a property into $\R^n$ with
the same Bellman function. This process is essentially trivial when we
use the Haar system in $\R^n$ introduced in Section 2, and it allows to do the ``induction in scales argument" at once, instead of once per
each $j=1,...,2^n-1$, which then introduces a dimensional constant of order $2^n$ in the estimates. We present this as the lemma named \emph{Good Bellman function Lemma}.
\begin{lemma}[Good Bellman Function Lemma]\label{GBFL}
Suppose there is a Bellman function, $B(X,Y)$ defined on domain $\mathfrak{D}$, that has the size property:
\begin{equation}0\leq B(X,Y)\leq A\,b(X,Y)\,,\label{GSize}\end{equation}
and the convexity property:
\begin{equation}B(X,Y)-\frac{B(X^1,Y^1)+B(X^2,Y^2)}{2}\geq C(X,Y,X^1,Y^1,X^2,Y^2,M)\,,\label{GConvex}\end{equation}
for $(X,Y)$ and $(X^{1,2},Y^{1,2})\in\mathfrak{D}\,,$
where $2X=X^1+X^2,$ and $Y=\frac{Y^1+Y^2}{2}+M\,.$
Furthermore, suppose the given Bellman function prove the certain dyadic property in $\R\,,$ that is,
for all dyadic interval $I\in\mathcal{D}\,,$
\begin{equation}\sum_{J\in\mathcal{D}(I)} C(X_J,Y_J,X_J^1,Y_J^1,X_J^2,Y_J^2,M_J)\leq |\,I|\,A\,b(X_I,Y_I)\,.
\label{DyadPro1}\end{equation}
Then, the extended property (\ref{DyadPro1}) to $\R^n$, that is
\begin{align}
 \sum_{Q\in\mathcal{D}^n(Q')}\sum_{j:E_{j,Q}\subseteq E_{i,Q`}}&C(X_{E_{j,Q}},Y_{E_{j,Q}},
X_{E^1_{j,Q}},Y_{E^1_{j,Q}},X_{E^2_{j,Q}},Y_{E^2_{j,Q}},M_{E_{j,Q}})\nonumber\\
&\leq |E_{i,Q'}|\,A\,b(X_{E_{i,Q'}},Y_{E_{i,Q'}})\label{DyadPron}\,,
\end{align}
is provided by the same Bellman function, $B\,,$ and checking $(X_{E_{j,Q}},Y_{E_{j,Q}})$ and $(X_{E^{1,2}_{j,Q}},Y_{E^{1,2}_{j,Q}})$ belong to $\mathfrak{D}\,,$
where $$\big(X_{E_{j,Q}},Y_{E_{j,Q}}\big)=\bigg(\frac{X_{E^1_{j,Q}}+X_{E^2_{j,Q}}}{2},M_{E_{i,Q}}+\frac{Y_{E^1_{j,Q}}+Y_{E^2_{j,Q}}}{2}\bigg)\,.$$
\end{lemma}
Note that the variables in (\ref{GConvex}), $X,~X^{1,2},~Y,~Y^{1,2},$ and $M$ can be considered as arbitrary tuples. The number of tuples depend on the given Bellman function.
\begin{proof} Without loss of generality we assume that $i=1\,.$
Let us assume that $(X_{E_{j,Q}},Y_{E_{j,Q}}),$ and $(X_{E^{1,2}_{j,Q}},Y_{E^{1,2}_{j,Q}})$ are on the domain $\mathfrak{D}\,,$
where $$\big(X_{E_{j,Q}},Y_{E_{j,Q}}\big)=\bigg(\frac{X_{E^1_{j,Q}}+X_{E^2_{j,Q}}}{2},M_{E_{j,Q}}+\frac{Y_{E^1_{j,Q}}+Y_{E^2_{j,Q}}}{2}\bigg)\,,$$
for all $Q\in\mathcal{D}^n(Q')$ and $j=1,...,2^n-1\,.$ Then, using the size condition (\ref{GSize}) and the convexity condition (\ref{GConvex}) we have,
for fixed dyadic cube $Q'$ and $i$,
\begin{align*}
A|E_{1,Q'}|&b(X_{E_{1,Q'}},Y_{E_{1,Q'}})\geq |E_{1,Q'}|B(X_{E_{1,Q'}},Y_{E_{1,Q'}})\\
&\geq \frac{|E_{1,Q'}|}{2}\sum_{l=1}^2B(X_{E^l_{1,Q'}},Y_{E^l_{1,Q'}})\\
&\q+C(X_{E_{1,Q'}},Y_{E_{1,Q'}},
X_{E^1_{1,Q'}},Y_{E^1_{1,Q'}},X_{E^2_{1,Q'}},Y_{E^2_{1,Q'}},M_{E_{1,Q'}})\\
&=\sum_{j=2}^{3}|E_{j,Q'}|B(X_{E_{j,Q'}},Y_{E_{j,Q'}})\\
&\q+C(X_{E_{j,Q'}},Y_{E_{j,Q'}},
X_{E^1_{i,Q'}},Y_{E^1_{i,Q'}},X_{E^2_{i,Q'}},Y_{E^2_{i,Q'}},M_{E_{i,Q'}})\,.\end{align*}
If we iterate this process $n-1$ times more, we get:
\begin{align*}
 A|Q'|&b(X_{Q'},Y_{Q'})\\
&\geq \sum_{j=2^{n-1}}^{2^n-1}\sum_{l=1}^2|E^l_{j,Q'}|B(X_{E^l_{j,Q'}},Y_{E^l_{j,Q'}})\\
&\q+\sum_{j=1}^{2^n-1} C(X_{E_{j,Q'}},Y_{E_{j,Q'}},X_{E^1_{j,Q'}},Y_{E^1_{j,Q'}},X_{E^2_{j,Q'}},Y_{E^2_{j,Q'}},M_{E_{j,Q'}})\,,
\end{align*}
Due to our construction of Haar system, for all $j=2^{n-1},2^{n-1}+1,...,2^n-1,$ and $l=1,2,$ $E^l_{j,Q'}$'s are mutually disjointed and
$|E_{j,Q'}^l|=|Q'|/2^n$ i.e. $\{E^1_{j,Q'},E^2_{j,Q'}\}_{j=2^n-1}^{2^{n-1}}$ is the set $\mathcal{D}^n_1(Q')$ of dyadic sub-cubes of $Q'$. Thus,
\begin{align*}
 A|Q'|&b(X_{Q'},Y_{Q'})\geq |Q'|B(X_{Q'},Y_{Q'})\\
&\geq \sum_{k=1}^{2^n}|Q'_k|B(X_{Q'_k},Y_{Q'_k})\\
&\q+\sum_{j=1}^{2^n-1} C(X_{E_{j,Q'}},Y_{E_{j,Q'}},X_{E^1_{j,Q'}},Y_{E^1_{j,Q'}},X_{E^2_{j,Q'}},Y_{E^2_{j,Q'}},M_{E_{j,Q'}})\,,
\end{align*}
where $Q'_k$'s are enumerations of $2^n$ dyadic sub-cubes of $Q'\,.$ Iterating this procedure and using the fact $B\geq 0$ yields that
\begin{align*}\sum_{Q\in\mathcal{D}^n(Q')}&\sum_{j=1}^{2^n-1}C(X_{E_{j,Q}},Y_{E_{j,Q}},X_{E^1_{j,Q}},Y_{E^1_{j,Q}},X_{E^2_{j,Q}},Y_{E^2_{j,Q}},M_{E_{j,Q}})\\
&\leq A|Q'|b(X_{Q'},Y_{E_{Q'}})\,,\end{align*} which completes the proof.
\end{proof}

We now state several multivariable versions of Embedding theorems and weight inequalities, that we will need to prove our main Theorems.
One can find the proof of Lemma \ref{lmwce} in \cite{NTV}.

\begin{lemma}\label{lmwce}
The following function
$$B(F,f,u,Y)=4A\bigg(F-\frac{f^2}{u+Y}\bigg)$$ is defined on domain
$\mathfrak{D}$ which is given by
$$\mathfrak{D}=\big\{(F,f,u,Y)\in\R^4\big|\,F,f,u,Y>0\textrm{ and }
f^2\leq Fu,~Y\leq u\big\}\,,$$ and $B$ satisfies the following
size and convexity property in $\mathfrak{D}$: $
0\leq B(F,f,u,Y)\leq 4A\,F\,,$
and for all $(F,f,u,Y),\,(F_1,f_1,u_1,Y_1)$ and
$(F_2,f_2,$ $u_2,Y_2)\in\mathfrak{D}\,,$
\begin{equation*}
B(F,f,u,Y)-\frac{B(F_1,f_1,u_1,Y_1)+B(F_2,f_2,u_2,Y_2)}{2}
\geq\frac{f^2}{u^2}M\,,
\end{equation*}
where
$$(F,f,u,Y)=\bigg(\frac{F_1+F_2}{2},\frac{f_1+f_2}{2},\frac{u_1+u_2}{2},M+\frac{Y_1+Y_2}{2}\bigg)\textrm{ and } m\geq 0\,.$$
\end{lemma}
Replacing
\begin{align*}X_{E_{j,Q'}}&=\big(\langle f^2\rangle_{E_{j,Q'}},\langle fw^{1/2}\rangle_{E_{j,Q'}},\langle w\rangle_{E_{j,Q}}\big)\,,\\
Y_{E_{j,Q'}}&=\frac{1}{|E_{j,Q'}|}\sum_{Q\in\mathcal{D}^n(Q')}\sum_{i:E_{i,Q}\subseteq E_{j,Q'}}\alpha_{i,Q}u^2_{j,Q}\,,
\end{align*}and
\begin{equation*}
M_{E_{j,Q'}}=\frac{1}{|E_{j,Q'}|}\alpha_{j,Q'}u^2_{j,Q'}\,,
\end{equation*}
in Lemma \ref{GBFL} and using Lemma \ref{lmwce}, we have the following Theorem whose one-dimensional version can be found in \cite{NTV}.
\begin{theorem}[Multivariable Version of Weighted Carleson Embedding Theorem]\label{MWCE}
Let $w$ be a weight and $\big\{\alpha_{j,Q}\big\}_{Q,j},$ $Q\in\mathcal{D}^n,~j=1,...,2^n-1\,,$ be a
sequence of nonnegative numbers such that for all dyadic cubes
$Q'\in\mathcal{D}^n$  and a positive constant
$A>0$,
\begin{equation}
\frac{1}{\big|E_{i,Q'}\big|}\sum_{Q\in\mathcal{D}^n(Q')}\sum_{\,j:E_{j,Q}\subseteq
E_{i,Q'}}\alpha_{j,Q}\langle w\rangle^2_{E_{j,Q}}\leq A\langle
w\rangle_{E_{i,Q'}}\,.\label{MBE:e1}\end{equation} Then for all positive $f\in L^2_{\R^n}$ \begin{equation}
\sum_{Q\in\mathcal{D}^n}\sum_{j=1}^{2^n-1}\alpha_{j,Q}\langle
fw^{1/2}\rangle^2_{E_{j,Q}}\leq CA\|f\|^2_{L^2_{\R^n}}\label{MBE:e2}
\end{equation}
holds with some constant $C>0\,.$
\end{theorem}

Similarly, one can prove the following theorem with Lemma \ref{GBFL} and the Bellman function appeared in \cite{Pet2}.

\begin{theorem}[Multivariable Version of Petermichl's the Bilinear Embedding Theorem]\label{PMBE}
Let $w$ and $v$ be weights so that $\langle w\rangle_{Q'}\langle w\rangle_{Q'}<A$ and $\big\{\alpha_{j,Q}\big\}_{Q,j}$ be
a sequence of nonnegative numbers such that, for all dyadic cubes
$Q'\in\mathcal{D}^n$ and $i=1,...,2^n-1\,,$ the three inequalities below holds
with some constant $A>0\,,$
\begin{align*}
&\frac{1}{\big|E_{i,Q'}\big|}\sum_{Q\in\mathcal{D}^n(Q')}\sum_{\,j:E_{j,Q}\subseteq E_{i,Q'}}\frac{\alpha_{j,Q}}{\langle w\rangle_{E_{j,Q}}}\leq A\langle v\rangle_{E_{i,Q'}} \\
&\frac{1}{\big|E_{i,Q'}\big|}\sum_{Q\in\mathcal{D}^n(Q')}\sum_{\,j:E_{j,Q}\subseteq E_{i,Q'}}\frac{\alpha_{j,Q}}{\langle v\rangle_{E_{j,Q}}}\leq A\langle w\rangle_{E_{i,Q'}}\\
&\frac{1}{\big|E_{i,Q'}\big|}\sum_{Q\in\mathcal{D}^n(Q')}\sum_{\,j:E_{j,Q}\subseteq
E_{i,Q'}}\alpha_{j,Q}\leq A\,.
\end{align*}
Then for all $f\in L^2_{\R^n}(w)$ and $g\in L^2_{\R^n}(v)$
$$\sum_{Q\in\mathcal{D}^n(Q')}\sum_{j=1}^{2^n-1}\alpha_{j,Q}\langle f\rangle_{E_{j,Q},w}\langle g\rangle_{E_{j,Q},v}\leq CA\|f\|_{L^2_{\R^n}(w)}\|g\|_{L^2_{\R^n}(v)}$$
holds with some constant $C>0\,.$
\end{theorem}

Changing $\alpha_{j,Q}$, $f$ and $g$ by $\alpha_{j,Q}\langle v\rangle_{E_{j,Q}}\langle w\rangle_{E_{j,Q}}|E_{j,Q}|$, $fw^{-1/2}$ and $gv^{-1/2}$ respectively in Theorem \ref{PMBE}, we can get the following Corollary.
\begin{corollary}[Multivariable Version of the Bilinear Embedding Theorem]\label{MBE}
Let $w$ and $v$ be weights so that $\langle w\rangle_{Q'}\langle w\rangle_{Q'}<A$ and $\big\{\alpha_{j,Q}\big\}_{Q,j}$ be
a sequence of nonnegative numbers such that, for all dyadic cubes
$Q'\in\mathcal{D}^n$ and $i=1,...,2^n-1\,,$ the three inequalities below holds
with some constant $A>0\,,$
\begin{align*}
&\frac{1}{\big|E_{i,Q'}\big|}\sum_{Q\in\mathcal{D}^n(Q')}\sum_{\,j:E_{j,Q}\subseteq E_{i,Q'}}\alpha_{j,Q}\langle v\rangle_{E_{j,Q}}|{E_{j,Q}}|\leq A\langle v\rangle_{E_{i,Q'}} \\
&\frac{1}{\big|E_{i,Q'}\big|}\sum_{Q\in\mathcal{D}^n(Q')}\sum_{\,j:E_{j,Q}\subseteq E_{i,Q'}}\alpha_{j,Q}\langle w\rangle_{E_{j,Q}}|{E_{j,Q}}|\leq A\langle w\rangle_{E_{i,Q'}}\\
&\frac{1}{\big|E_{i,Q'}\big|}\sum_{Q\in\mathcal{D}^n(Q')}\sum_{\,j:E_{j,Q}\subseteq
E_{i,Q'}}\alpha_{j,Q}\langle w\rangle_{E_{j,Q}}\langle
v\rangle_{E_{j,Q}}|{E_{j,Q}}|\leq A\,.
\end{align*}
Then for all $f,g\in L^2_{\R^n}$
$$\sum_{Q\in\mathcal{D}^n(Q')}\sum_{j=1}^{2^n-1}\alpha_{j,Q}\langle fw^{1/2}\rangle_{E_{j,Q}}\langle gv^{1/2}\rangle_{E_{j,Q}}|{E_{j,Q}}|\leq CA\|f\|_{L^2_{\R^n}}\|g\|_{L^2_{\R^n}}$$
holds with some constant $C>0\,.$
\end{corollary}

We now state several propositions include the multidimensional analogues to the corresponding one-dimensional results in \cite{Be} and \cite{Pet2} for both regular and anisotropic cases.

One can find the 1-dimensional analogue of the following proposition and the associated Bellman function lemma in \cite{Be}. Repeating the proof of Lemma \ref{GBFL} with $X_{j,Q'}=\big(\langle w\rangle_{E_{j,Q'}},\langle w^{-1}\rangle_{E_{j,Q'}}\big),$
$Y_{E_{j,Q'}}=\frac{1}{A|E_{j,Q'}|}\sum_{Q\in\mathcal{D}^n(Q')}\sum_{i:E_{i,Q}\subseteq E_{j,Q'}}\alpha_{i,Q},$ and $M_{j,Q'}=\frac{\alpha_{j,Q'}}{A}\,,$
and corresponding Bellman function lemma will return the following proposition.

\begin{proposition}\label{wp1}
Let $w$ be a weight, so that $w^{-1}$ is also a weight. Let $\alpha_{j,Q}$ be a
Carleson sequence of nonnegative numbers i.e., there is a constant
$A>0$ such that, for all $Q'\in\mathcal{D}^n$ and
$i=1,...,2^n-1\,,$
\begin{equation}\frac{1}{\big|E_{i,Q'}\big|}\sum_{Q\in\mathcal{D}^n(Q')}\sum_{j:E_{j,Q}\subseteq
E_{i,Q'}}\alpha_{j,Q}\leq A\,.\label{wp1:ee}\end{equation} Then,
for all $Q'\in\mathcal{D}^n$ and $i=1,...,2^n-1\,,$
\begin{equation}\frac{1}{\big|E_{i,Q'}\big|}\sum_{Q\in\mathcal{D}^n(Q')}\sum_{\,j:E_{j,Q}\subseteq E_{i,Q'}}\frac{\alpha_{j,Q}}
{\langle w^{-1}\rangle_{E_{j,Q}}}\leq 4A\langle
w\rangle_{E_{i,Q'}}\,,\label{wp1:e1}\end{equation} and if $w\in
A^d_2$ then for any $Q'\in\mathcal{D}^n$ and $i=1,...,2^n-1\,,$ we
have
\begin{equation}
\frac{1}{|E_{i,Q'}|}\sum_{Q\in\mathcal{D}^n(Q')}\sum_{\,j:E_{j,Q}\subseteq
E_{i,Q'}}\langle w\rangle_{E_{j,Q}}\alpha_{j,Q}\leq 4\cdot
2^{2(n-1)}A[w]_{A^d_2}\langle w\rangle_{E_{i,Q'}}\,.\label{wp1:e2}
\end{equation}
Furthermore, if $w\in A^R_2$ then for any $Q'\in\mathcal{D}^n$ and
$i=1,...,2^{n-1}\,,$ we have
\begin{align}
\frac{1}{|E_{i,Q'}|}\sum_{Q\in\mathcal{D}^n(Q')}\sum_{\,j:E_{j,Q}\subseteq
E_{i,Q'}}\langle w\rangle_{E_{j,Q}}\alpha_{j,Q}\leq
4A[w]_{A^R_2}\langle w\rangle_{E_{i,Q'}}\,.\label{dwp1:e2}
\end{align}
\end{proposition}
Observe that in the case $w\in A^R_2$ then $$\frac{1}{\langle w^{-1}\rangle_{E_{j,Q}}}\geq \frac{\langle w\rangle_{E_{j,Q}}}{[w]_{A_2^R}}\,.$$ Now (\ref{dwp1:e2}) follows from (\ref{wp1:e1}). Observe if $w \in A_2^d$ then
\begin{align*}
[w]_{A_2}&\geq \langle w\rangle_Q\langle w^{-1}\rangle_Q\\
&\geq \bigg(\frac{|E_{j,Q}|}{|Q|}\bigg)^2\langle w\rangle_{E_{j,Q}}\langle w^{-1}\rangle_{E_{j,Q}}=2^{-2(n-1)}\langle w\rangle_{E_{j,Q}}\langle w^{-1}\rangle_{E_{j,Q}}\,.
\end{align*}
Thus, we can have (\ref{wp1:e2}) from (\ref{dwp1:e2}).
We refer to \cite{Pet2} for the corresponding one-dimensional result and the associated Bellman function lemma to Proposition \ref{wp2}. Lemma \ref{GBFL} with $X_{E_{j,Q'}}=\big(\langle w\rangle_{E_{j,Q'}},\langle w^{-1}\rangle_{E_{j,Q}}\big)$, $Y_{E_{j,Q'}}
=M_{E_{j,Q'}}=0\,$ and the associated Bellman function lemma yield Proposition \ref{wp2}.
\begin{proposition}\label{wp2}
There exist a positive constant $C$ so that for all weight $w$ and
$w^{-1}$ and for all dyadic cubes $Q'\in\mathcal{D}^n$ and
$i=1,...,2^n-1$:
\begin{equation}\frac{1}{\big|E_{i,Q'}\big|}\sum_{Q\in\mathcal{D}^n(Q')}\sum_{\,j:E_{j,Q}\subseteq E_{i,Q'}}\frac{\big(\langle
w\rangle_{E^1_{j,Q}}-\langle w\rangle_{E^2_{j,Q}}\big)^2}{\langle
w\rangle^3_{E_{j,Q}}}\big|E_{j,Q}\big|\leq C\langle
w^{-1}\rangle_{E_{i,Q'}}\label{wp2:m}\end{equation} and, if $w\in A^d_2\,,$ the following
inequality holds for all dyadic cubes $Q'\in\mathcal{D}^n\,$ and
$i=1,...,2^n-1:$
\begin{align}&\frac{1}{\big|E_{i,Q'}\big|}\sum_{Q\in\mathcal{D}^n(Q')}\sum_{\,j:E_{j,Q}\subseteq E_{i,Q'}}\Bigg(\frac{\langle
w\rangle_{E^1_{j,Q}}-\langle w\rangle_{E^2_{j,Q}}}{\langle
w\rangle_{E_{j,Q}}}\Bigg)^2\big|E_{j,Q}\big|\langle
w^{-1}\rangle_{E_{j,Q}}\nonumber\\
&\qquad\leq C2^{2(n-1)}[w]_{A^d_2}\langle
w^{-1}\rangle_{E_{i,Q'}}\,.\label{wp2:e1}\end{align} Moreover, if
$w\in A^R_2\,,$ the following inequality holds for all dyadic
cubes $Q'\in\mathcal{D}^n\,$ and $i=1,...,2^n-1:$
\begin{align}&\frac{1}{\big|E_{i,Q'}\big|}\sum_{Q\in\mathcal{D}^n(Q')}\sum_{\,j:E_{j,Q}\subseteq E_{i,Q'}}\Bigg(\frac{\langle
w\rangle_{E^1_{j,Q}}-\langle w\rangle_{E^2_{j,Q}}}{\langle
w\rangle_{E_{j,Q}}}\Bigg)^2\big|E_{j,Q}\big|\langle
w^{-1}\rangle_{E_{j,Q}}\nonumber\\
&\qquad\leq C[w]_{A^R_2}\langle
w^{-1}\rangle_{E_{i,Q'}}\,.\label{dwp2:e1}\end{align}
\end{proposition}

The similar observations in Proposition \ref{wp1} yields (\ref{wp2:e1}) and (\ref{dwp2:e1}).
The following generalizes the result that appeared in \cite{Be} to
the multidimensional regular and anisotropic cases. With same changes in Proposition \ref{wp2} and associated Bellman function lemma, one can prove the following.

\begin{proposition}\label{wp3}
There exist a positive constant $C$ so that for all weight $w$ and
$w^{-1}$ and for all dyadic cubes $Q'\in\mathcal{D}^n$ and
$i=1,...,2^n-1:$
\begin{align}
&\frac{1}{\big|E_{i,Q'}\big|}\sum_{Q\in\mathcal{D}^n(Q')}\sum_{\,j:E_{j,Q}\subseteq
E_{i,Q'}}\Bigg(\frac{\langle w\rangle_{E^1_{j,Q}}-\langle
w\rangle_{E^2_{j,Q}}}{\langle
w\rangle_{E_{j,Q}}}\Bigg)^2\big|E_{j,Q}\big|\langle
w\rangle^{1/4}_{E_{j,Q}}\langle
w^{-1}\rangle^{1/4}_{E_{j,Q}}\nonumber\\
&\qquad\leq C\langle w\rangle^{1/4}_{E_{i,Q'}}\langle
w^{-1}\rangle^{1/4}_{E_{i,Q'}}\label{wp3:m}\end{align} and, if $w\in A^d_2\,,$
the following inequality holds for all dyadic cubes
$Q'\in\mathcal{D}^n\,:$
\begin{align}&\frac{1}{\big|E_{i,Q'}\big|}\sum_{Q\in\mathcal{D}^n(Q')}\sum_{\,j:E_{j,Q}\subseteq E_{i,Q'}}\Bigg(\frac{\langle
w\rangle_{E^1_{j,Q}}-\langle w\rangle_{E^2_{j,Q}}}{\langle
w\rangle_{E_{j,Q}}}\Bigg)^2\big|E_{j,Q}\big|\langle
w\rangle_{E_{j,Q}}\langle w^{-1}\rangle_{E_{j,Q}}\nonumber\\
&\qquad\leq C2^{2(n-1)}[w]_{A^d_2}\,.\label{wp3:e1}\end{align} Moreover,
if $w\in A^R_2\,,$ the following inequality holds for all dyadic
cubes $Q'\in\mathcal{D}^n\,:$
\begin{align}&\frac{1}{\big|E_{i,Q'}\big|}\sum_{Q\in\mathcal{D}^n(Q')}\sum_{\,j:E_{j,Q}\subseteq E_{i,Q'}}\Bigg(\frac{\langle
w\rangle_{E^1_{j,Q}}-\langle w\rangle_{E^2_{j,Q}}}{\langle
w\rangle_{E_{j,Q}}}\Bigg)^2\big|E_{j,Q}\big|\langle
w\rangle_{E_{j,Q}}\langle w^{-1}\rangle_{E_{j,Q}}\nonumber\\
&\qquad\leq
C[w]_{A^R_2}\,.\label{dwp3:e1}\end{align}
\end{proposition}

The single variable version of the
following proposition first appeared in \cite{Wi1}. In \cite{Per},
one can also find a Bellman function proof of a similar result
which can be extended to the doubling measure case.

\begin{proposition}[Wittwer's sharp version of Buckley's
inequality]\label{wp4} There exist a positive constant $C$ so that
for all weight $w\in A^d_2$ and all dyadic cubes
$Q'\in\mathcal{D}^n$ and $i=1,...,2^n-1:$
\begin{align}&\frac{1}{\big|E_{i,Q'}\big|}\sum_{Q\in\mathcal{D}^n(Q')}\sum_{\,j:E_{j,Q}\subseteq E_{i,Q'}}\Bigg(\frac{\langle
w\rangle_{E^1_{j,Q}}-\langle w\rangle_{E^2_{j,Q}}}{\langle
w\rangle_{E_{j,Q}}}\Bigg)^2\big|E_{j,Q}\big|\langle
w\rangle_{E_{j,Q}}\nonumber\\
&\qquad\leq C 2^{2(n-1)}[w]_{A^d_2}\langle
w\rangle_{E_{i,Q'}}\,,\label{wp4:e1}\end{align} and for all
weight $w\in A^R_2$ and all dyadic cubes $Q'\in\mathcal{D}^n$ and
$i=1,...,2^n-1:$
\begin{align}&\frac{1}{\big|E_{i,Q'}\big|}\sum_{Q\in\mathcal{D}^n(Q')}\sum_{\,j:E_{j,Q}\subseteq E_{i,Q'}}\Bigg(\frac{\langle
w\rangle_{E^1_{j,Q}}-\langle w\rangle_{E^2_{j,Q}}}{\langle
w\rangle_{E_{j,Q}}}\Bigg)^2\big|E_{j,Q}\big|\langle
w\rangle_{E_{j,Q}}\nonumber\\
&\qquad\leq C[w]_{A^R_2}\langle
w\rangle_{E_{i,Q'}}\,.\label{dwp4:e1}\end{align}
\end{proposition}
The same choice of $X_{E_{j,Q'}},Y_{E_{j,Q'}},$ and $M_{E_{j,Q'}}$ with Proposition \ref{wp4} prove (\ref{wp4:e1}).
The inequality (\ref{dwp4:e1}) can be seen by using the domain
$\mathfrak{D}=\big\{(u,v)\in\R^2\,|\,u,v>0\textrm{ and }1\leq
uv\leq [w]_{A^R_2}\big\}\,.$

\section{Proof of Theorem \ref{MA}}
We are going to prove Theorem \ref{MA} only when $p=2\,,$ and following the one-dimensional proof discovered by Beznosova \cite{Be}. The
sharp extrapolation theorem \cite{DGPerPet} returns immediately
the other cases $(1<p<\infty)$. For the case $p=2$ we use the
duality arguments. Precisely, it is sufficient to prove the
inequality
\begin{equation}
\langle \pi_b(fw^{-1/2}),gw^{1/2}\rangle \leq
C(n)[w]_{A^d_2}\|b\|_{BMO^d_{\R^n}}\|f\|_{L^2_{\R^n}}\|g\|_{L^2_{\R^n}}\,.\label{PT:e1}
\end{equation}
\begin{proof}
Using the orthogonal Haar system (\ref{PN:e5}), we can split the
left hand side of (\ref{PT:e1}) as follows.
\begin{align}
\langle \pi_b(f&w^{-1/2}),gw^{1/2}\rangle=\sum_{Q\in\mathcal{D}^n}\sum_{j=1}^{2^n-1}\langle b,h_{j,Q}\rangle\langle fw^{-1/2}\rangle_{E_{j,Q}}\langle gw^{1/2},h_{j,Q}\rangle\nonumber\\
=&\sum_{Q\in\mathcal{D}^n}\sum_{j=1}^{2^n-1}\langle b,h_{j,Q}\rangle\langle fw^{-1/2}\rangle_{E_{j,Q}}\langle gw^{1/2},H^w_{j,Q}\rangle\frac{1}{\sqrt{|E_{j,Q}|}}\label{PT:e2}\\
&~~+\sum_{Q\in\mathcal{D}^n}\sum_{j=1}^{2^n-1}\langle b,h_{j,Q}\rangle\langle
fw^{-1/2}\rangle_{E_{j,Q}}\langle
gw^{1/2},A^w_{j,Q}\chi_{E_{j,Q}}\rangle\frac{1}{\sqrt{|E_{j,Q}|}}\,\,\,.\label{PT:e3}
\end{align}
We are going to prove that both sum (\ref{PT:e2}) and (\ref{PT:e3}) are
bounded with a bound that depends linearly on both $[w]_{A^d_2}$ and $\|b\|_{BMO^d_{\R^n}}\,.$ Using Cauchy-Schwarz
inequality, for the term (\ref{PT:e2}), we have
\begin{align}
&\sum_{Q\in\mathcal{D}^n}\sum_{j=1}^{2^n-1}\,\langle b,h_{j,Q}\rangle\langle fw^{-1/2}\rangle_{E_{j,Q}}\langle gw^{1/2},H^w_{j,Q}\rangle\frac{1}{\sqrt{|E_{j,Q}|}}\nonumber\\
&\q\leq\Bigg(\sum_{Q\in\mathcal{D}^n}\sum_{j=1}^{2^n-1}\frac{\langle gw^{1/2},H^w_{j,Q}\rangle^2}{|E_{j,Q}|\langle w\rangle_{E_{j,Q}}}\Bigg)^{1/2}\nonumber\\
&\hspace{3cm}\times\Bigg(\sum_{Q\in\mathcal{D}^n}\sum_{j=1}^{2^n-1}\langle b,h_{j,Q}\rangle^2\langle fw^{-1/2}\rangle^2_{E_{j,Q}}\langle w\rangle_{E_{j,Q}}\Bigg)^{1/2}\nonumber\\
&\q\leq\,
\|g\|_{L^2_{\R^n}}\Bigg(\sum_{Q\in\mathcal{D}^n}\sum_{j=1}^{2^n-1}\langle
b,h_{j,Q}\rangle^2\langle fw^{-1/2}\rangle^2_{E_{j,Q}}\langle
w\rangle_{E_{j,Q}}\Bigg)^{1/2}\label{PT:e6}\,.
\end{align}
Here the inequality (\ref{PT:e6}) follows from (\ref{PN:e7}). We
now claim that the sum in (\ref{PT:e6}) is bounded by
$C[w]^2_{A^d_2}\|b\|^2_{BMO^d_{\R^n}}\|f\|^2_{L^2_{\R^n}}\,,$
which will be provided by the Multivariable Version of the
Weighted Carleson Embedding Theorem \ref{MWCE}. with the embedding
condition: For all $Q'\in\mathcal{D}^n,$
\begin{align}
&\frac{1}{\big|E_{i,Q'}\big|}\sum_{Q\in\mathcal{D}^n(E_{i,Q'})}\sum_{j:E_{j,Q}\subseteq
E_{i,Q}}\langle w^{-1}\rangle^2_{E_{j,Q}}\langle w\rangle
_{E_{j,Q}}\langle b, h_{j,Q}\rangle^2\nonumber\\
&\qquad\leq C[w]^2_{A^d_2}\|b\|^2_{BMO^d_{\R^n}}\langle
w^{-1}\rangle_{E_{i,Q'}}\,.\label{PT:e7}
\end{align}
Since, for all $Q\in\mathcal{D}^n\,,$ $2^{2(n-1)}[w]_{A^d_2}\geq
\langle w\rangle_{E_{j,Q}}\langle w^{-1}\rangle_{{E_{j,Q}}}\,,$
The Embedding condition (\ref{PT:e7}) can be seen as follows.
\begin{align}
&\sum_{Q\in\mathcal{D}^n(E_{i,Q'})}\sum_{j:E_{j,Q}\subseteq E_{i,Q}}\,\langle w^{-1}\rangle^2_{{E_{j,Q}}}\langle w\rangle_{E_{j,Q}}\langle b, h_{j,Q}\rangle^2\nonumber\\
&\qquad\leq\,2^{2(n-1)}[w]_{A^d_2}\sum_{Q\in\mathcal{D}^n(Q')}\sum_{\,j:E_{j,Q}\subseteq E_{i,Q'}}\langle w^{-1}\rangle_{E_{j,Q}}\langle b, h_{j,Q}\rangle^2\nonumber\\
&\qquad\leq 4\cdot
2^{7(n-1)}[w]^2_{A^d_2}\|b\|^2_{BMO^d_{\R^n}}w^{-1}(E_{i,Q'})\label{PT:e8}\,.
\end{align}
Here the inequality (\ref{PT:e8}) follows from (\ref{BMO:e2}) and
Proposition \ref{wp1} applied to $\alpha_{j,Q}=\langle
b,h_{j,Q}\rangle^2$, $A=3^{n-1}\|b\|^2_{BMO^d_{\R^n}}$
and $v=w^{-1}\,.$ This estimates finishes the estimate of the term (\ref{PT:e2}) with $C\approx 2^{7(n-1)/2}\,.$

We now turn to estimate of the term (\ref{PT:e3}). In order
to estimate the term (\ref{PT:e3}), we need to show that
\begin{align}
&\sum_{Q\in\mathcal{D}^n}\sum_{j=1}^{2^n-1}\langle b,h_{j,Q}\rangle\langle
fw^{-1/2}\rangle_{E_{j,Q}}\langle
gw^{1/2}\rangle_{E_{j,Q}}\,A^w_{j,Q}\sqrt{|E_{j,Q}|}\nonumber\\
&\qquad\leq
C[w]_{A^d_2}\|b\|_{BMO^d_{\R^n}}\|f\|_{L^2_{\R^n}}\|g\|_{L^2_{\R^n}}\,,\label{PT:e9}
\end{align}
and this is provided by the following three embedding conditions
due to the Multivariable Version of the Bilinear Embedding Theorem
\ref{MBE}: For all $Q'\in\mathcal{D}^n$ and $i=1,...,2^n-1\,,$
\begin{align}
\frac{1}{\big|E_{i,Q'}\big|}&\sum_{Q\in\mathcal{D}^n(Q')}\sum_{\,j:E_{j,Q}\subseteq E_{i,Q'}}\big|\langle b,h_{j,Q}\rangle\,A^w_{j,Q}\big|\sqrt{|E_{j,Q}|}\langle w\rangle_{E_{j,Q}}\langle w^{-1}\rangle_{E_{j,Q}}\nonumber\\
&\leq C(n)[w]_{A^d_2}\|b\|_{BMO^d_{\R^n}}\,,\label{PT:e10}\\
\frac{1}{\big|E_{i,Q'}\big|}&\sum_{Q\in\mathcal{D}^n(Q')}\sum_{\,j:E_{j,Q}\subseteq E_{i,Q'}}\big|\langle b,h_{j,Q}\rangle\,A^w_{j,Q}\big|\sqrt{|E_{j,Q}|}\langle w\rangle_{E_{j,Q}}\nonumber\\
&\leq C(n)[w]_{A^d_2}\|b\|_{BMO^d_{\R^n}}\langle w\rangle_{E_{i,Q'}}\,,\label{PT:e11}\\
\frac{1}{\big|E_{i,Q'}\big|}&\sum_{Q\in\mathcal{D}^n(Q')}\sum_{\,j:E_{j,Q}\subseteq
E_{i,Q'}}\big|\langle
b,h_{j,Q}\rangle\,A^w_{j,Q}\big|\sqrt{|E_{j,Q}|}\langle
w^{-1}\rangle_{E_{j,Q}}\nonumber\\
&\leq
C(n)[w]_{A^d_2}\|b\|_{BMO^d_{\R^n}}\langle
w^{-1}\rangle_{E_{i,Q'}}\,.\label{PT:e12}
\end{align}
Proposition \ref{wp3} makes it easy to prove the embedding
condition (\ref{PT:e10}). Using Cauchy-Schwarz inequality,
\begin{align}
\sum_{Q\in\mathcal{D}^n(Q')}&\sum_{\,j:E_{j,Q}\subseteq E_{i,Q'}}\,\big|\langle b,h_{j,Q}\rangle\,A^w_{j,Q}\big|\sqrt{|E_{j,Q}|}\langle w\rangle_{E_{j,Q}}\langle w^{-1}\rangle_{E_{j,Q}}\\
&\leq\Bigg(\sum_{Q\in\mathcal{D}^n(Q')}\sum_{\,j:E_{j,Q}\subseteq E_{i,Q'}}\langle b,h_{j,Q}\rangle^2\langle w\rangle_{E_{j,Q}}\langle w^{-1}\rangle_{E_{j,Q}}\Bigg)^{1/2}\\
&\hspace{0.5cm} \times \Bigg(\sum_{Q\in\mathcal{D}^n(Q')}\sum_{\,j:E_{j,Q}\subseteq E_{i,Q'}} \big(A^w_{j,Q}\big)^2|E_{j,Q}|\langle w\rangle_{E_{j,Q}}\langle w^{-1}\rangle_{E_{j,Q}}\Bigg)^{1/2}\nonumber\\
&\leq\,2^{n-1}[w]^{1/2}_{A^d_2}\Bigg(\sum_{Q\in\mathcal{D}^n(Q')}\sum_{\,j:E_{j,Q}\subseteq E_{i,Q'}}\langle b,h_{j,Q}\rangle^2\Bigg)^{1/2}\\
&\hspace{0.5cm} \times \Bigg(\sum_{Q\in\mathcal{D}^n(Q')}\sum_{\,j:E_{j,Q}\subseteq E_{i,Q'}} \big(A^w_{j,Q}\big)^2|E_{j,Q}|\langle w\rangle_{E_{j,Q}}\langle w^{-1}\rangle_{E_{j,Q}}\Bigg)^{1/2}\nonumber\\
&\leq\,C\,2^{2(n-1)}[w]_{A^d_2}\big|E_{i,Q'}\big|^{1/2}\Bigg(\sum_{Q\in\mathcal{D}^n(Q')}\sum_{\,j:E_{j,Q}\subseteq E_{i,Q'}}\langle b,h_{j,Q}\rangle^2\Bigg)^{1/2}\label{PT:e14}\\
&\leq\,C\,2^{5(n-1)/2}[w]_{A^d_2}\|b\|_{BMO^d_{\R^n}}\big|E_{i,Q'}\big|\label{PT:e15}\,.
\end{align}
Here we use (\ref{wp3:e1}) for the inequality (\ref{PT:e14}) and
the fact that $b\in BMO^d_{\R^n}$ for the inequality
(\ref{PT:e15}). We also use Cauchy-Schwarz inequality for the
inequality (\ref{PT:e11}), then
\begin{align}
&\sum_{Q\in\mathcal{D}^n(Q')}\sum_{\,j:E_{j,Q}\subseteq E_{i,Q'}}\big|\langle b,h_{j,Q}\rangle\,A^w_{j,Q}\big|\sqrt{|E_{j,Q}|}\langle w\rangle_{E_{j,Q}}\nonumber\\
&\leq\Bigg(\sum_{Q\in\mathcal{D}^n(Q')}\sum_{\,j:E_{j,Q}\subseteq E_{i,Q'}}\langle b,h_{j,Q}\rangle^2\langle w\rangle_{E_{j,Q}}\Bigg)^{1/2}\nonumber\\
&\hspace{3cm} \times \Bigg(\sum_{Q\in\mathcal{D}^n(Q')}\sum_{\,j:E_{j,Q}\subseteq E_{i,Q'}} \big(A^w_{j,Q}\big)^2|E_{j,Q}|\langle w\rangle_{E_{j,Q}}\Bigg)^{1/2}\nonumber\\
&\leq C2^{5(n-1)/2}\|b\|_{BMO^d_{\R^n}}[w]^{1/2}_{A_2^d}\langle w\rangle_{E_{i,Q'}}^{1/2}\nonumber\\
&\hspace{3cm}\times\Bigg(\sum_{Q\in\mathcal{D}^n(Q')}\sum_{\,j:E_{j,Q}\subseteq E_{i,Q'}} \big(A^w_{j,Q}\big)^2|E_{j,Q}|\langle w\rangle_{E_{j,Q}}\Bigg)^{1/2}\label{PT:e16}\\
&\leq C\,2^{7(n-1)/2}\|b\|_{BMO^d_{\R^n}}[w]_{A_2^d}\langle
w\rangle_{E_{i,Q'}}\label{PT:e17}\,.
\end{align}
Inequality (\ref{PT:e16}) and (\ref{PT:e17}) follow by
(\ref{wp1:e2}) and Proposition \ref{wp4} respectively. Similarly, we can
establish inequality (\ref{PT:e12}) with Proposition \ref{wp2}. To sum up, we can establish the inequality (\ref{PT:e1}) with a
constant $C(n)\approx 2^{7(n-1)/2}\,.$ Furthermore, if we replace
$[w]_{A_2^d}$ by $[w]_{A^R_2}$ and $\|b\|_{BMO^d}$ by
$\|b\|_{BMO^R}$ then we can establish proof of the dimension free
estimate in Theorem \ref{MA}.
\end{proof}

\section{Final Remarks}
\begin{remark}
The martingale transforms in $\R$ is defined by
$$T_{\sigma}f:=\sum_{I\in\mathcal{D}([0,1])}\sigma_I\langle f,h_I\rangle h_I,$$
where $\sigma_I=\pm 1$. It is also known to be a dyadic analog of singular integral operators. In \cite{Wi1}, the author presented a linear estimate
of the martingale transform on the weighted Lebesgue space $L^2(w)$, that is
$$\|T_{\sigma}f\|_{L^2(w)}\leq C[w]_{A_2}\|f\|_{L^2(w)}\,,$$ for $w\in A_2$ and $f\in L^2(w)\,.$ With the Wilson's Haar system, we can also define the multivariable martingale transform:
\begin{equation}\label{MMT}
 T_{\sigma}f=\sum_{Q\in\mathcal{D}^n([0,1]^n)}\sum_{j=1}^{2^n-1}\sigma_{E_{j,Q}}\langle f,h_{j,Q}\rangle h_{j,Q}\,,
\end{equation}
where $\sigma_{E_{j,Q}}$ assumes the values $\pm 1$ only. It was also considered in \cite{DPetV} to search for the $L^p$
estimate of the Beurling-Ahlfors operator.

 In fact, we already present all the tools to extend the result of \cite{Wi1} to (\ref{MMT}). In \cite{Wi1}, the one dimensional analogue of the inequality:
\begin{align}
 \frac{1}{|E_{i,Q'}|}&\sum_{Q\in\mathcal{D}^n(Q')}\sum_{j:E_{j,Q}\subseteq E_{i,Q'}} \big(\langle w^{-1}\rangle_{E^1_{j,Q}}-\langle w\rangle_{E^2_{j,Q}}\big)^2
|E_{j,Q}|\langle w\rangle_{E_{j,Q}}\nonumber\\
&\leq C(n)[w]_{A_2}^2\langle w^{-1}\rangle_{E_{i,Q'}}\label{mmte1}
\end{align}
was proved by using the sharp estimate of the dyadic square function in $L^2(w)$. Alternatively, one can also have the inequality (\ref{mmte1}) using
(\ref{wp1:e2}), (\ref{wp3:e1}) and $[w]_{A_2}=[w^{-1}]_{A_2}.$ By the Cauchy-Schwartz inequality and Proposition \ref{wp3} we have the following inequality:
\begin{align}
\sum_{Q\in\mathcal{D}^n(Q')}&\sum_{j:E_{j,Q}\subseteq E_{i,Q'}}
(\langle w\rangle_{E^1_{j,Q}}-\langle w\rangle_{E^2_{j,Q}})(\langle w^{-1}\rangle_{E^1_{j,Q}}-\langle w^{-1}\rangle_{E^2_{j,Q}})|E_{j,Q}|\nonumber\\
\leq& C(n)[w]_{A_2}|E_{i,Q'}|\,.\label{mmte2}
\end{align}
One can find the Bellman function proof of the single variable version of the following inequality.
\begin{align}
 &\sum_{Q\in\mathcal{D}^n(Q')}\sum_{j:E_{j,Q}\subseteq E_{i,Q'}}\frac{(\langle w\rangle_{E^1_{j,Q}}-\langle w\rangle_{E^2_{j,Q}})
(\langle w^{-1}\rangle_{E^1_{j,Q}}-\langle w^{-1}\rangle_{E^2_{j,Q}})}{\langle w^{-1}\rangle_{E_{j,Q}}}|E_{j,Q}|\nonumber\\
&\qquad\quad\leq C(n)[w]_{A_2}w(E_{i,Q'})\,.\label{mmte3}
\end{align}
Using Lemma \ref{GBFL} we can have the inequality (\ref{mmte3}). However, more simply we have
\begin{align}
 &\sum_{Q\in\mathcal{D}^n(Q')}\sum_{j:E_{j,Q}\subseteq E_{i,Q'}}\frac{(\langle w\rangle_{E^1_{j,Q}}-\langle w\rangle_{E^2_{j,Q}})
(\langle w^{-1}\rangle_{E^1_{j,Q}}-\langle w^{-1}\rangle_{E^2_{j,Q}})}{\langle w^{-1}\rangle_{E_{j,Q}}}|E_{j,Q}|\nonumber\\
&\q\leq \bigg(\sum_{Q\in\mathcal{D}^n(Q')}\sum_{j:E_{j,Q}\subseteq E_{i,Q'}}(\langle w\rangle_{E^1_{j,Q}}-\langle w\rangle_{E^2_{j,Q}})^2|E_{j,Q}|
\langle w\rangle_{E_{j,Q}}\bigg)^{1/2}\label{mmte4}\\
&\quad\times\bigg(\sum_{Q\in\mathcal{D}^n(Q')}\sum_{j:E_{j,Q}\subseteq E_{i,Q'}}\bigg(\frac{\langle w^{-1}\rangle_{E^1_{j,Q}}-\langle w^{-1}\rangle_{E^2_{j,Q}}}{\langle
w^{-1}\rangle_{E_{j,Q}}}\bigg)^2|E_{j,Q}|\langle w\rangle_{E_{j,Q}}\bigg)^{1/2}\label{mmte5}\,.
\end{align}
Using the inequality (\ref{wp4:e1}) for the term (\ref{mmte4}) and (\ref{wp2:e1}) interchanging the role of $w$ and $w^{-1}$ for the term (\ref{mmte5}), we
can get the inequality (\ref{mmte3}). By repeating the proof presented in the Section 5 in \cite{Wi1} with previous observations and the sharp extrapolation theorem \cite{DGPP}, we have that for $1<p<\infty$ there exists constants $C(n,p)$ only depending on $p$ and dimension $n$ and $C$ which doesn't depend on the dimensional constant such
that
\begin{equation}
 \|T_{\sigma}\|_{L^p_{\R^n}(w)\ra L^p_{\R^n}(w)}\leq C(n,p)[w]_{A^d_p}^{\max\{1,\frac{1}{p-1}\}}\,,
\end{equation}
for all weights $w\in A^d_p$ and
\begin{equation}
 \|T_{\sigma}\|_{L^2_{\R^n}(w)\ra L^2_{\R^n}(w)}\leq C[w]_{A^R_2}\,,
\end{equation}
for all weights $w\in A^R_2\,.$
\end{remark}

\begin{remark}
Most recently, in \cite{CrMP}, it has been showed that the norm of the dyadic square function,
\begin{equation}S_df(x):=\bigg(\sum_{Q\in\mathcal{D}^n}\big(\langle f\rangle_Q-\langle f\rangle_{\widehat{Q}}\big)^2\chi_Q(x)\bigg)^{1/2}\,,\label{dsf}\end{equation} where $\widehat{Q}$ is the dyadic parent of $Q$, is bounded by $[w]_{A_p}^{\max\{\frac{1}{2},\frac{1}{p-1}\}}$ and the exponent is the best possible. The dyadic square function (\ref{dsf}) can be written in terms of Wilson's Haar system (see \cite{Wilt}), that is
\begin{equation}
S_df(x):=\Bigg(\sum_{Q\in\mathcal{D}^n}\frac{1}{|Q|}\bigg(\sum_{j=1}^{2^n-1}\langle f,h_{j,Q}\rangle^2\bigg)\chi_{Q}(x)\Bigg)^{1/2}\,.
\end{equation}
\end{remark}

\begin{remark}
Let us consider the difference between $\pi^s_b$ and $\pi_b$ as an operator,
\begin{equation}
\pi_b^sf(x)-\pi_bf(x)=\sum_{Q\in\mathcal{D}^n}\sum_{j=1}^{2^n-1}\big(\langle f\rangle_Q-\langle f\rangle_{E_{j,Q}}\big)\langle b,h_{j,Q}\rangle h_{j,Q}(x)\,.\label{dpsp}
\end{equation}
For fixed $Q$ and $j>1$, $\langle f\rangle_Q-\langle f\rangle_{E_{j,Q}}$ can be written by
$$\sum_{i:E_{i,Q}\supsetneq E_{j,Q}}\langle f,h_{i,Q}\rangle h_{i,Q}(E_{j,Q})\,.$$
Thus, the difference operator (\ref{dpsp}) can be estimated by
\begin{equation}
\big|\pi_b^sf(x)-\pi_bf(x)\big|\leq C(n)\sum_{Q\in\mathcal{D}^n}\sum_{j=1}^{2^n-1}\sum_{i:E_{i,Q}\supsetneq E_{j,Q}}|\langle f,h_{i,Q}\rangle|\,|\langle b,h_{j,Q}\rangle|\frac{\chi_{E_{j,Q}}(x)}{|E_{j,Q}|}\,.\label{fsoap}
\end{equation}
The right hand side of the inequality (\ref{fsoap}) can be considered as a finite sum of the compositions of adjoint of paraproducts and dyadic shifting operators, that is denoted by $\pi_b^{\ast}H_{\tau}\,.$ It is known that the dyadic shifting operators obey linear bounds in $L^2(w)$. We have shown that the paraproduct ($\pi_b$) obey linear bounds in $L^2(w)$, so the adjoint of paraproduct also does. One can easily expect that its composition, $\pi^{\ast}_bH_{\tau}$ obey quadratic bounds in $L^2(w)\,.$ However, $\pi^{\ast}_bH_{\tau}$ obey linear bounds in $L^2(w)$ \cite{Ch} but $H_{\tau}\pi_b^{\ast}$ doesn't. See \cite{Ch} for more detailed arguments in the one dimensional case.
Then the difference operator $\pi_b^s-\pi_b$ obeys the linear bound in $L^2_{\R^n}(w)$. Furthermore, this observation and the linear bound for the paraproduct associated to the Wilson's Haar basis ($\pi_b$) yield the linear bound for the paraproduct associated to the standard tensor product Haar basis ($\pi_b^s$).
\end{remark}

\noindent \emph{Acknowledgment.} This work is part of the author's Ph.D dissertation \cite{Ch}. The author like to thank his graduate advisor Mar\'{i}a Cristina Pereyra for her suggestions and
helpful interaction.

\bibliographystyle{amsplain}

\end{document}